\documentclass[10pt,a4paper]{amsart}

\usepackage{amsfonts,amsmath,amssymb,graphics,setspace}

\input xy
\xyoption{all}

\newtheorem{thm}{Theorem}[section]
\newtheorem{lemma}[thm]{Lemma}
\newtheorem{cor}[thm]{Corollary}
\newtheorem{pro}[thm]{Proposition}
\newtheorem{example}[thm]{Example}
\newtheorem{remark}[thm]{Remark}
\newtheorem{ddef}[thm]{Definition}

\def\limind{\setbox1=\hbox{\oalign{\vadjust{\vskip -2pt}%
     \rm lim\cr \vadjust{\vskip -2pt}
      \hidewidth$-\mkern -12mu\rightarrow$\hidewidth\cr}}
       \hbox{\box1}}

\def\ext{\mathop{\rm Ext}\nolimits}

\def\tor{\mathop{\rm Tor}\nolimits}
\def\hom{\mathop{\rm Hom}\nolimits}
\def\shhom{\mathop{\mathcal Hom}\nolimits}
\def\spec{\mathop{\rm Spec}}

\def\depth{\mathop{\rm depth}\nolimits}
\def\cd{\mathop{\rm cd}\nolimits}
\def\reg{\mathop{\rm reg}\nolimits}

\def\a{{\alpha}}

\def\g{{\gamma}}

\def\ZZ{{\bf Z}}
\def\NN{{\bf N}}

\def\fin{\mathop{\rm end}}
\def\ara{\mathop{\rm ara}}

\def\ann{\mathop{\rm ann}\nolimits}
\def\coker{\mathop{\rm coker}}

\def\ass{\mathop{\rm Ass}\nolimits}
\def\Min{\mathop{\rm Min}\nolimits}
\def\height{\mathop{\rm ht}\nolimits}

\def\depth{\mathop{\rm depth}\nolimits}

\def\supp{\mathop{\rm Supp}}

\def\G{{\Gamma}}
\def\GG{{\bf \Gamma}}
\def\RG{{R\Gamma}}
\def\RGG{{R{\bf\Gamma}}}

\def\A{{\mathcal A}}

\def\C{{\mathcal C}}
\def\F{{\mathcal F}}
\def\E{{\mathcal E}}
\def\OO{{\mathcal O}}

\def\om{{\omega}}

\def\a{{\alpha}}

\def\g{{\gamma}}

\def\pp{{S_+}}
\def\ip{{\mathfrak p}}
\def\iq{{\mathfrak q}}
\def\Ip{{\mathfrak P}}
\def\Iq{{\mathfrak Q}}
\def\ia{{\mathfrak a}}

\def\im{{\mathfrak m}}
\def\idn{{\mathfrak n}}

\def\ra{{\rightarrow}}
\def\lra{{\longrightarrow}}
\def\fini{{$\quad\quad\Box$}}

\newcommand{\bd}{\begin{ddef}}
\newcommand{\ed}{\end{ddef}}
\newcommand{\bt}{\begin{thm}}
\newcommand{\et}{\end{thm}}
\newcommand{\bl}{\begin{lemma}}
\newcommand{\el}{\end{lemma}}
\newcommand{\bco}{\begin{cor}}
\newcommand{\eco}{\end{cor}}
\newcommand{\bp}{\begin{pro}}
\newcommand{\ep}{\end{pro}}
\newcommand{\bex}{\begin{example}}
\newcommand{\eex}{\end{example}}
\newcommand{\brm}{\begin{remark}}
\newcommand{\erm}{\end{remark}}
\newcommand{\bconj}{\begin{conj}}
\newcommand{\econj}{\end{conj}}

\newcommand{\ult}{\underline{T}}

\newcommand{\Ext}{\operatorname{Ext}}

\newcommand{\beqn}{\begin{eqnarray*}}
\newcommand{\eeqn}{\end{eqnarray*}}
\newcommand{\beq}{\begin{eqnarray}}
\newcommand{\eeq}{\end{eqnarray}}
\newcommand{\been}{\begin{enumerate}}
\newcommand{\eeen}{\end{enumerate}}

\begin{document}

\author{Marc Chardin}
\address{Institut de Math\'ematiques de Jussieu\\
UPMC, Boite 247, 4, place Jussieu, F-75252 Paris Cedex, France}
\email{chardin@math.jussieu.fr}

\author{Jean-Pierre Jouanolou}
\address{IRMA, Universit\'e de Strasbourg\\
7, rue Ren\'e Descartes
F-67084 Strasbourg, France}
\email{jean-pierre.jouanolou@math.unistra.fr}

\author{Ahad Rahimi}
\address{Department of Mathematics, Razi University,
Kermanshah, Iran
\&  School of Mathematics, IPM, P. O. Box 19395-5746, Tehran, Iran.}
\email{ahad.rahimi@razi.ac.ir}

\title[]
{The eventual stability of depth, associated primes and cohomology of a graded module}

\maketitle

\section*{Introduction}

The asymptotic stability of several homological invariants of the graded pieces 
of a graded module has attracted quite a lot of attention over the last decades. An early
important result was the proof by Brodmann of the eventual stabilization of the associated 
primes of the powers of an ideal in a Noetherian ring (\cite{Br1}). 

We provide in this text several stability results together with estimates of the degree 
from which it stabilizes. One of our initial goals was to obtain a simple proof of the
tameness result of Brodmann in \cite{Br2} for graded components of cohomology over rings of
dimension at most two. This is achieved in the last section, and gives a slight generalization 
of what is known, as our result (Theorem \ref{tamedim2}) applies to Noetherian rings of dimension at 
most two that are either local or the epimorphic image of a Gorenstein ring. Recall that Cutkosky and 
Herzog provided examples in \cite{CH} showing that tameness does not hold over rings of dimension three (even
over nice local such rings). 

Besides this result, we establish, for a graded module $M$ over a polynomial ring $S$
(in finitely many variables, with its standard grading)  over a 
commutative ring $R$, stability results for the depth and cohomological dimension of
graded pieces with
respect to a finitely generated $R$-ideal $I$.  It follows from our results that 
the cohomological dimension of $M_\mu$ with respect to $I$ is constant for $\mu >\reg (M)$,
and the depth with respect to $I$ is at least equal to its eventual value for $\mu >\reg (M)$ and
stabilizes when it reaches this value for some $M_\mu$ with $\mu >\reg (M)$. See \ref{assdepth} and
\ref{stabcd} for more precise results.

Recall that $\reg (M)\in {\bf Z}$ when $M\not= 0$ is finitely generated and
$R$ is Noetherian.

When $R$ is Noetherian,
$\ip\in\spec (R)$ is associated to $M_\mu$ for some $\mu$ if and only if $\ip =\Ip\cap R$ for
$\Ip$ associated to $M$ in $S$ and the set associated primes of $M_\mu$ is non decreasing 
for $\mu >\reg (M)$. It implies that this set eventually stabilizes when $M$ is finitely generated.

Before we establish these regularity results in sections 4 and 5, we prove several facts about
depth and cohomological
dimension with respect to a finitely generated ideal and about Castelnuovo-Mumford regularity of a 
graded module. Our definition of depth agrees with the one introduced by Northcott. These results are 
stated in a quite general setting and self contained proofs are given.
Our arguments are often at least as simple as the ones proposed under stronger hypotheses in classical 
references. 
We are in particular careful 
about separating statements  where a finiteness hypothesis is needed (notably in terms of finite generation, 
finite presentation, or Noetherianity) from others that do not require it.  
We show that 
several basic results on regularity hold without any finiteness hypothesis, and that many results 
on the asymptotic behaviour hold for modules of finite regularity. 

In Section 6, we give pretty general duality statements that encapsulate the Herzog-Rahimi spectral sequence
we use in the last section to derive tameness from our previous stability results. 

\section{Local cohomology and depth}

Let $A$ be a commutative ring (with unit) and $M$ a $A$-module. 
If $a=(a_1,\ldots ,a_r)$ is a $r$-tuple of elements of $A$, $K^\bullet (a;M)$ is the Koszul complex and $H^i(a;M)$ its $i$-th cohomology module. Also, $\C^\bullet_a (M)$ is the \v Cech complex. This complex is isomorphic to  $\limind_n K^\bullet (a_1^n,\ldots ,a_r^n;M)$. If $a$ and $b$  generates two ideals with same radical, then $H^i(\C^\bullet_a (M))\simeq H^i(\C^\bullet_b (M))$ for
all $i$. Moreover this isomorphism is graded (of degree 0) if $A$, $M$ and the ideals generated by $a$ and $b$ 
are graded. This for instance follows from \cite[1.2.3 and 1.4.1]{EGAIII}. It can also be proved in an elementary way as follows :
first notice that it is sufficient to prove that if $y\in \sqrt{(x_1,\ldots x_t)}$ then $H^i(\C^\bullet_{(x_1,\ldots x_t)}(M))\simeq H^i(\C^\bullet_{(x_1,\ldots x_t,y)} (M))$,
second show that $\C^\bullet_{(x_1,\ldots x_t)} (M_{y})$ is acyclic if $y\in \sqrt{(x_1,\ldots x_t)}$, and conclude using 
that $ \C^\bullet_{(x_1,\ldots x_t,y)} (M)$ is the mapping cone of the natural map
$\C^\bullet_{(x_1,\ldots x_t)} (M)\ra \C^{\bullet}_{(x_1,\ldots x_t)} (M_{y})$.

We will denote by $H^i_I (M)$ the $i$-th homology module of $H^i(\C^\bullet_a (M))$, if
$a$ generates the ideal $I$.

The $i$-th right derived functor of the left exact functor $H^0_I$ coincides with the functor $T^i(\hbox{---}):=\limind_n \ext^i_A(A/I^n,\hbox{---})$.
It coincides with  $H^i_I $ if and only if $H^i_I(N)=0$ whenever $i>0$ and $N$ is injective, and this holds if $A$
is Noetherian or $I$ is generated by a regular sequence.

If $X:=\spec (A)$ and $Y:=V(I)\subset X$, one has an isomorphism 
$$
H^i_I (M) \simeq H^i_Y (X,\widetilde{M}).
$$
 Indeed Serre affineness theorem  and Cartan-Leray theorem (see e.g. \cite{Ke} or \cite{Su}, and \cite[5.9.1]{Go}) provide isomorphisms 
\begin{equation}\label{SCL}
H^i_Y (X,\widetilde{M})\simeq H^i(\C^\bullet_{(a_1,\ldots a_r)} (M)) \simeq H^i(M\otimes_A^{\bf L} \C^\bullet_{(a_1,\ldots a_r)} (A)),
\end{equation}
as $\C^\bullet_{(a_1,\ldots a_r)} (A)$ is a complex of flat modules. These isomorphisms show that the functor $M\mapsto H^i_I (M)$ commutes with direct sums and filtered inductive limits and provide a spectral sequence
\begin{equation}\label{SSTor}
E^{p,q}_{2}=\tor_{-p}^{A}(M,H^q_I (A))\, \Rightarrow H^{p+q}_I (M).
\end{equation}

Also notice that the isomorphism $H^i_I (M)\simeq \limind_n H^i (a_1^n,\ldots ,a_r^n;M)$ shows that any element of $H^i_I (M)$ is
annihilated by a power of the ideal $I$.

\bd
If $I$ is a finitely generated $A$-ideal and $M$ an $A$-module, we set 
$$
\depth_I (M):=\max
\{ p\in \NN \cup \{ +\infty \} \ \vert \ H^i_I (M)=0, \forall i<p\} ,
$$
and
$$
\cd_I (M):=\max
\{ p\in \NN \cup \{ -\infty \} \ \vert \ H^p_I (M)\not= 0\} .
$$
\ed

In case there might be an ambiguity on the ring over which $I$ and $M$ are considered, 
we will use the notations $\depth_I^A (M)$ and $\cd_I^A (M)$.

Notice that, for any $A$-module $M$, $\cd_I (M)$ is bounded above by the minimal number
of generators of any ideal $J$ such that $\sqrt{J}=\sqrt{I}$ (this number is called the arithmetic rank 
of $I$ in $A$, $\ara_A (I)$).

\bl\label{KosCech}
If $I$ is generated by $a=(a_1,\ldots ,a_r)$,
$$
\depth_I (M)=\max 
\{ p\in \NN \cup \infty \ \vert \ H^i(a;M)=0, \forall i<p\} .
$$
\el

{\it Proof.} Let $d:=\max
\{ p\ \vert \ H^i(a;M)=0, \forall i<p\}$. Recall that for positive integers 
$l_i$, $H^i(a_1^{l_1},\ldots ,a_r^{l_r};M)=0$ if and only 
if $H^i(a;M)=0$. It follows that $H^i_I (M)=\limind_n H^i(a_1^n,\ldots ,a_r^n;M)=0$ if 
$H^i(a;M)=0$. Notice that $d=\infty$ if and only if
$d>r$, in which case $H^i(a;M)=H^i_I (M)=0$ for all $i$. Hence $\depth_I (M)\geq d$.  
We now assume $d<\infty$. As $I$ annihilates 
$H^i(a;M)$ for any $i$,
the totalisation of the complex $\C^\bullet_I K^\bullet (a;M)$ has cohomology isomorphic
to the one of $K^\bullet (a;M)$. It provides a spectral sequence
$$
{'E}_1^{p,q}=H^q_{I}K^p(a;M) \ \Rightarrow H^{p+q}(a;M).
$$
As $H^q_I K^p(a;M)=0$ for $q<d$, this in turn provides a natural into map $H^{d}(a;M)\ra H^d_{I}(M)$ which
shows that $\depth_I (M)\leq d$.\fini

\bco
If $I$ is a finitely generated $A$-ideal, then for any $A$-module $M$,
$$
\depth_I (M) =\min_{\ip \in V(I)}\{ \depth_{I_\ip} (M_\ip )\} .
$$
\eco

To show that this notion agrees with the one introduced by Northcott, we first prove a lemma.

\bl\label{depthNzd}
Let $N$ be a $A$-module and $a\in I$ a non zero divisor on $N$. Then
$$
\depth_I (N/aN)=\depth_I (N)-1.
$$
\el

{\it Proof.} Consider the exact sequence
$$
\xymatrix{
0\ar[r]&N\ar[r]^{\times a}&N\ar[r]&N/aN\ar[r]&0\\}
$$
and the induced long exact sequence on cohomology with support in $I$,
$$
\xymatrix{
\cdots\ar[r]&H^{i-1}_I (N)\ar[r]&H^{i-1}_I (N/aN)\ar[r]&H^{i}_I (N)\ar[r]^{\times a}&H^{i}_I (N)\ar[r]&\cdots\\}
$$
 and let $r:=\depth_I (N)$. The above sequence shows that $\depth_I (N/aN)\geq r-1$.
 Furthermore, if $r<\infty$, $H^{r-1}_I (N/aN)=0$ if and only if the multiplication by $a$ is injective on $H^{r}_I (N)$. But this does not hold since any element of $H^{r}_I (N)$ is annihilated by a power of $a$ and $H^{r}_I (N)\not= 0$ by definition.\fini

We will also use a version of the Dedekind-Mertens Lemma, that we now recall in its
general form, together with immediate corollaries that are useful in this text.

\bt\label{LDM} (Generalized Dedekind-Mertens Lemma, \cite[3.2.1]{Jo}). 
Let $A$ be a ring, $M$ be a $A$-module and $\ult$ a set of variables. For
$P\in A[\ult ]$ and $Q \in M[\ult ]$, let $c(P)$ be the $A$-ideal generated by the coefficients of $P$, 
$c(Q)$ be the submodule of $M$ generated by the coefficients of $Q$ and $\ell (Q)$ be the number of non zero coefficients of $Q$. 

Then one has the equality
$$
c(P)^{\ell (Q)-1}c(PQ)=c(P)^{\ell (Q)}c(Q).
$$
In particular, the kernel of the multiplication by $P$ in $M[\ult ]$ is supported in $V(c(P))$.
\et

\bco\label{DMA}
Let $A$ be a ring, $I=(a_0,\ldots ,a_p)$ be a $A$-ideal and $M$ be a $A$-module. Set
$\xi :=a_0+a_1 T+\cdots +a_p T^p\in A[T]$. Then
$$
\ker (\!\!\xymatrix{M[T]\ar^{\times \xi}[r]& M[T]\\}\!\! )\subset H^0_I (M[T])=H^0_I (M)[T].
$$
\eco

Let $S=R[X_1,\ldots  ,X_n]$ be a polynomial ring over a commutative ring $R$ and
set $\pp:=(X_1,\ldots  ,X_n)$.

\bco\label{kglf} 
Let $M$ be a graded $S$-module. Set $\ell :=T_1X_1+\cdots +T_nX_n$  with $\deg T_i=0$. 
Then the kernel of the map,
$$
\xymatrix{M[T_1,\ldots ,T_n ]\ar^{\times \ell}[r]&M[T_1,\ldots ,T_n](1)\\}
$$
is a graded $S[T_1,\ldots ,T_n ]$-submodule of $H^0_\pp (M)[T_1,\ldots ,T_n ]$.
\eco

\bco\label{glt} Consider indeterminates $(U_{i,j})_{1\leq i,j\leq n}$,  $\xi_i :=\sum_{1\leq j\leq n}U_{i,j}X_j$,
$\Delta :=\det (U_{i,j})_{1\leq i,j\leq n}$ and $R':=R[(U_{i,j})_{1\leq i,j\leq n}]_\Delta$. Let 
$S':=R'[X_1,\ldots  ,X_n]$ and set $M':=M\otimes_R R'$ for any $S$-module $M$.
Then $(\xi_1 ,\ldots ,\xi_n )$ is $M'$-regular
  off $V(\pp ')=V(\xi_1 ,\ldots ,\xi_n )$. 
\eco

The following proposition shows that the above definition of depth agrees with the one introduced by Northcott
in \cite{No}.

\bp\label{NorthDepth}
Let $r\geq 1$ be an integer and $I$ be a finitely generated $A$-ideal. The following are equivalent,

(1) $\depth_I (M)\geq r$,

(2) There exists a faithfully flat extension $B$ of $A$ and a regular sequence 
$f_1,\ldots ,f_r$ on $B\otimes_A M$ contained in $IB$,

(3) There exists a polynomial extension $B$ of $A$ and a regular sequence 
$f_1,\ldots ,f_r$ on $B\otimes_A M$ contained in $IB$,

(4) There exists  a regular sequence 
$f_1,\ldots ,f_r$ on $M[T_1,\ldots ,T_r]$, where the $T_i$'s are variables, contained in $IA[T_1,\ldots ,T_r]$.
\ep

{\it Proof.} The implications ${\rm (4)}\Rightarrow {\rm (3)}\Rightarrow {\rm (2)}$ are 
trivial. Furthermore (2) implies that $H^i_{IB} (B\otimes_A M )=0$ for $i<r$ using Lemma \ref{depthNzd}, which in turn
implies (1) since $H^i_{IB} (B\otimes_A M )\simeq B\otimes_A H^i_{I} (M )$ because $B$ is flat 
over $A$.

Finally (1) implies (4) by  induction on $r$, using  Lemma  \ref{depthNzd} and Corollary \ref{DMA}.\fini

\brm\label{CorNorthDepth}
Let $r\geq 1$ be an integer and $I$ be a finitely generated $A$-ideal. 
If $\depth_I (M)\geq r$ and $f_1,\ldots ,f_s\in IB$ is a regular sequence on
$B\otimes_A M$ for some flat extension $B$ of $A$, then $s\leq r$ and there exists a faithfully flat extension $C$ of
$B$ and $f_{s+1},\ldots ,f_r\in IC$ such that $f_1,\ldots ,f_r$ is regular on $C\otimes_A M$.
\erm

\section{Castelnuovo-Mumford regularity}

Let $S$ be a finitely generated standard graded algebra over a commutative ring $R$.
Recall that for a graded $S$-module $M$
$$
\reg (M):=\sup\{ \mu\ \vert\ \exists i,\ H^i_{S_+}(M)_{\mu -i}\not= 0\} ,
$$
with the convention that $\sup \emptyset =-\infty$.

\bl\label{3.1}
Let $M$ be a graded $S$-module. Consider the following properties,

(i) $M_\mu =0$ for $\mu\gg 0$,

(ii) $M=H^0_{S_+}(M)$,

(iii) $H^i_{S_+}(M)=0$ for $i>0$.

Then $(i)\Rightarrow (ii)\Rightarrow (iii)$, $(ii)\Rightarrow (i)$ if  $M$ is finitely generated or $\reg (M)<\infty$, and $(iii)\Rightarrow (ii)$ if $M_\mu =0$
for $\mu \ll 0$.
\el

{\it Proof.} $(i)\Rightarrow (ii)$ is clear since any homogeneous element in $M$ is killed by a power of any
element of $S_+$ if $(i)$ holds. If $(ii)$ holds, then $M_x =0$ for any $x\in S_+$, hence the \v Cech complex on
generators of $S_1$ (as an $R$-module) is concentrated in homological degree $0$, which shows $(iii)$. 

If $(ii)$ holds,
any element in $M$ is killed by a power of $S_+$, hence if $M$ is finitely generated by $(m_t)_{t\in T}$,  any generator $m_t$
is killed by $S_+^{N_t}$, for some $N_t \in \NN$. It then follows that any element in $M$ of degree bigger than
$\max_{t\in T}\{ \deg (m_t )+N_t \} $ is $0$.  If $\reg (M)<\infty$, $(ii)\Rightarrow (i)$ follows trivially from the definition of $\reg (M)$. 

If $(iii)$ holds set $N:=M/H^0_{S_+} (M)$.
The exact sequence $0\ra H^0_{S_+} (M)\ra M\ra N\ra 0$ gives rise to a long exact sequence in local cohomology showing that $H^i_{S_+} (N)=0$ 
for all $i$. As $\depth_{S_+}(N)=+\infty$, Lemma \ref{KosCech} shows that $N=S_+N$. This implies that $N=0$ as $N_\mu =0$ for $\mu\ll 0$.
\fini

The following two propositions extend classical results on regularity. 

\bp
Let $\ell \geq 1$ and $m$ be integers and $M$ be a graded $S$-module. 

If $H^i_{S_+}(M)_{m -i}= 0$ for $i\geq \ell$,
then $H^i_{S_+}(M)_{\mu -i}= 0$ for $\mu \geq m$ and $i\geq \ell$.

Assume that $H^i_{S_+}(M)_{m -i}= 0$ for all $i$ and let $\mu \geq m$.
Then $H^i_{S_+}(M)_{\mu -i}= 0$ for $i>0$ and  $H^0_{S_+}(M)_{\mu}=M_\mu /S_1M_{\mu -1}$. 

\ep
 
{\it Proof.} We may assume that $S=R[X_1,\ldots ,X_n]$ is a polynomial ring over $R$. We then prove the assertion
by induction on $n$. 

When $n=0$,  $M=H^0_\pp (M)$ and $H^i_\pp (M)=0$ for $i\not= 0$, and the claim follows in both cases.

Next assume that $n\geq 1$ and the assertion is true for $n-1$ over any commutative ring. Let
$$
\xi :=X_1T^{n-1}+\cdots +X_{n-1}T+X_n \in R[T,X_1,\ldots ,X_n].
$$

The Dedekind-Mertens Lemma implies that 
$$
\ker (\xi :M[T](-1)\ra M[T])\subset H^0_\pp (M[T]),
$$
by  Corollary \ref{DMA}. But $H^i_\pp (N[T])=H^i_\pp (N)[T]$ for any
$S$-module $N$ and any $i$. Hence, replacing $R$ by $R[T]$, $M$ by $M[T]$ and $X_n$ by $\xi$, we may assume that 
$$
K:=\ker (X_n :M(-1)\ra M)\subset H^0_\pp (M).
$$
We then have $H^i_\pp (K)=0$ for $i\not= 0$ and the exact sequence
$$
\xymatrix{
0\ar[r]&K\ar[r]&M(-1)\ar^{\times X_n}[r]&M\ar[r]&Q\ar[r]&0\\
}
$$
induces for all $i$ an exact sequence
$$
\xymatrix{
H^i_\pp (M)(-1)\ar^{\times X_n}[r]&H^i_\pp (M)\ar[r]&H^i_\pp (Q)\ar[r]&H^{i+1}_\pp (M)(-1).\\
}
$$
For  $i+j=m$, the equalities
$$
H^i_\pp (M)_j =0\ {\rm and}\ H^{i+1}_\pp (M)(-1)_j=H^{i+1}_\pp (M)_{j-1}=0
$$
imply $H^i_\pp (Q)_j =0$. Hence $H^i_\pp (Q)_j =0$ for $i\geq \ell$ and $i+j=m$. As $Q$ is annihilated by $X_n$, setting $\idn :=(X_1,\ldots ,X_{n-1})$ one has
$H^i_\pp (Q)=H^i_\idn (Q)$ for all $i$. Applying the recursion hypothesis to the $R[X_1,\ldots ,X_{n-1}]$-module $Q$,
it follows that 
$$
H^i_\pp (Q)_j =H^i_\idn (Q)_j =0,\quad \forall i\geq \ell ,\ \forall j\geq m-i.
$$
Hence $X_n : H^i_\pp (M)_{j-1}\ra H^i_\pp (M)_{j}$ is onto for $i\geq \ell $ and $i+j\geq m$. As $H^i_\pp (M)_{m-i}=0$ for
$i\geq \ell$, this proves our claim. \fini


\brm
Notice that the exact sequence $0\ra S_+ M\ra M\ra M/S_+ M\ra 0$ induces an exact sequence
$$
0\ra H^0_{S_+}(S_+ M)\ra H^0_{S_+}(M)\ra M/S_+ M\ra H^1_{S_+}(S_+ M)\ra H^1_{S_+}(M)\ra 0
$$
\erm

\bp
If $S=R[X_1,\ldots ,X_n]$ is a polynomial ring, then for any graded $S$-module $M$

(i) $\reg (M)=\sup \{ \mu\ \vert\ \exists i,\ \tor_i^{S}(M,R)_{\mu +i}\not= 0\}$,

(ii) $\reg (M)=\sup \{ \mu\ \vert\ \exists i,\ H_i(X_1,\ldots ,X_n;M,)_{\mu +i}\not= 0\}$,

(iii) $\reg (M)=\sup \{ \mu\ \vert\ \exists j,\ H^j(X_1,\ldots ,X_n;M,)_{\mu -j}\not= 0\}$,

(iv)  $\reg (M)=\sup \{ \mu\ \vert\ \exists j,\ \ext_S^j(R,M,)_{\mu -j}\not= 0\}$.

\ep

In particular, $M$ is generated in degrees at most $\reg (M)$ (when $\reg (M)=-\infty$, it means that
 $M$ is generated in degrees at most $\mu$, for any $\mu\in {\bf Z}$).\medskip

{\it Proof.} We first show this equality if $\reg (M)<\infty$. Let $K_\bullet (M):=K_\bullet (X_1,\ldots
,X_n;M)$ and $K^\bullet (M):=K^\bullet (X_1,\ldots
,X_n;M)$. As $K_\bullet (M)=K_\bullet (S)\otimes_S M$,  $K_\bullet (M)=\hom_S (K_\bullet (S),M)$ and $K_\bullet (S)$ is a free $S$-resolution of $R$, it
follows that $H_i (K_\bullet (M))\simeq \tor_i^{S}(M,R)$ and $H^j (K^\bullet (M))\simeq \ext_S^{j}(R,M)$. 
Furthermore, the complexes $K^\bullet (M)$ and $K_\bullet (M)$ are isomorphic, up to a shift in homological degree and internal degree : $K^\bullet (M)\simeq K_{n-\bullet}(M)(n)$, proving that 
$$
\ext_S^{j}(R,M)_{\mu -j}\simeq H^j (K^\bullet (M))_{\mu -j}\simeq H_{n-j} (K_\bullet (M))_{\mu -j+n}\simeq\tor_{n-j}^{S}(M,R)_{\mu -j+n}
$$ 
and the equivalence of the four items.

To prove (ii), the double complex $\C_{(X_1,\ldots ,X_n)}^\bullet K_\bullet (M)$
gives rise to two spectral sequences whose first terms are respectively
$$
'E_1^{p,q}=K_{-p}(X_1,\ldots ,X_n;H^q_{S_+} (M)),\quad
'E_2^{p,q}=\tor_{-p}^{S}(H^q_{S_+}(M),R)\quad
$$
and
$$
''E_1^{p,q}=\C_{(X_1,\ldots ,X_n)}^p \tor_{-q}^{S}(M,R),\quad
''E_2^{p,q}=H^p_{S_+}(\tor_{-q}^{S}(M,R)).
$$
Recall that $K_\bullet (M)_{X_i}$ is acyclic for any $i$,
hence $''E_1^{p,q}=0$ for $p\not= 0$ which implies that 
$''E_\infty^{p,q}=0$ for $p\not= 0$ and $''E_\infty^{0,q}\simeq
\tor_{-q}^{S}(M,R)$.

On the other hand, $('E_1^{p,q})_{\mu}=0$ for $\mu >\reg (M)-q-p$
as $H^q_{S_+} (M)$ lives in degrees at most $\reg (M)-q$. It follows
first that $\tor_{-j}^{S}(M,R)$ lives in degrees at most $\reg (M)-j$ showing
that $\tor_i^{S}(M,R)_{\mu +i}=0$ for $\mu >\reg (M)$.

To conclude, choose $j$ such that $H^j_{S_+}(M)_{\reg (M)-j}\not= 0$.
Set $\mu :=\reg (M)-j+n$ and notice that $('E_1^{p,q})_\mu =0$ when
$p+q=-n+j+1$. As  $'E_1^{p,q}=0$ for $p<-n$ it follows
that $0\not= H^j_{S_+}(M)_{\reg (M)-j}=('E_1^{-n,j})_\mu \simeq ('E_\infty^{-n,j})_\mu$
which shows that $''E_\infty^{0,j-n}\simeq\tor_{n-j}^{S}(M,R)_{\reg (M)+n-j}\not= 0$.

To finish the proof, we must show that $\reg (M)<\infty$ if there exists $\mu_0$ such that
$\tor_i^{S}(M,R)_{\mu}= 0$ for all $i$ and $\mu \geq \mu_0$. 

We first show that in this case,
there exists a graded free $S$-resolution $F_\bullet$ of $M$ with $F_i=\oplus_{j\in I_i} S(-d_{ij})$,
and $d_{ij}<\mu_0$ for all $i$ and $j$. Notice that if $M$ is graded and $(M/S_{+}M)_{>\nu }=0$, then 
$M$ is generated in degree at most $\nu$, showing the existence of a graded epimorphism 
$\phi : F_0\ra M$ with $F_0$ as claimed.  The exact sequence,
$$
0\ra \ker (\phi )\ra F_0 \ra M\ra0
$$
gives rise to another
$$
\tor_1^{S}(M,R)\ra \ker (\phi )/S_{+}\ker (\phi) \ra F_0/S_{+}F_0 
$$
and proves the existence of $F_1$ as claimed such that $F_1\ra \ker (\phi )$ is a graded epimorphism.
As $\tor_j^{S}(M,R)\simeq \tor_{j-1}^{S}(\ker (\phi ),R)$, for $j\geq 2$, the conclusion follows by 
induction on $i$. 

Finally, it suffices to remark that if $F_\bullet$ is a graded resolution as above, then 
$H^i_{S_{+}}(M)_\mu \simeq H_{n-i}(H^n_{S_+} (F_\bullet )_\mu )$ vanishes in degree
bigger than $-n+\max_{i,j}\{ d_{ij}\}\leq -n+\mu_0$, as $H^n_{S_+} (F_\bullet )_\mu =0$ for
all $i$ for such a $\mu$. Hence $\reg (M)\leq \mu_0$.
\fini

\bl\label{localreg}
For any  graded $S$-module
$N$, 
$$
\reg (N)=\sup_{\ip \in \spec (R)}\{ \reg (N\otimes_{R}R_\ip )\}.
$$
Furthermore, $\reg (N)=\reg (N\otimes_{R}R_\ip )$ for some $\ip\in \spec (R)$ if $\reg (N)<\infty$.
\el

{\it Proof.} For $\ip\in \spec (R)$,
$$
H^i_{S_+}(N\otimes_{R}R_\ip )_\mu \simeq (H^i_{S_+}(N)_\mu )\otimes_{R}R_\ip 
$$
from which the claim directly follows.\fini

\section{Depth of the graded components of a graded module}

As in the previous section,  $S$ is a finitely generated standard graded algebra over a commutative ring $R$.

\medskip
\bp\label{assdepth}
Let $I$ be a finitely generated  $R$-ideal, $M$ be a graded $S$-module and $d$ be an integer. 
 Assume that $\depth_I (M_{\mu })\geq d$ for  $\mu\gg 0$. Then 
 
 (i) $\depth_I (M_\mu )\geq d$ for any $\mu$ such that $H^q_{S_+}(M)_\mu =0$ for $q<d$,
 
 (ii) if $\depth_I (M_\mu )=d$ for some $\mu$ such that $H^q_{S_+}(M)_\mu =0$ for $q\leq d$, then
 $\depth_I (M_\nu )\leq d$ for any $\nu\geq \mu$.
\ep

{\it Proof.} Let $a=(a_1,\ldots ,a_r)$ be generators of $I$. Recall that  $H^p(a;M)_\mu =H^p(a;M_\mu )$, as $I$ is a $R$-ideal. 
Hence, if $\depth_I (M_{\mu })\geq d$ for  $\mu\gg 0$, the $S$-modules 
$H^p(a;M)$ are supported in $V(S_{+})$ for $p<d$, by Lemma \ref{KosCech} and $(i)\Rightarrow (ii)$ in Lemma \ref{3.1}. Comparing the two 
spectral sequences obtained from the double complex $\C^{\bullet}_{S_+}K^{\bullet}(a;M)$,
computing the hypercohomology ${\bf H}^{\bullet}$ of $K^{\bullet}$ we obtain on one hand 
that ${\bf H}^i\simeq H^i(a;M)$ for $i<d$ and ${\bf H}^d\simeq H^0_{S_+}(H^d(a;M))$. On the other hand, one has a 
spectral sequence 
$$
{'E}_1^{p,q}=H^q_{S_+}K^p \ \Rightarrow {\bf H}^{p+q}
$$
which shows  that $({\bf H}^i)_\mu =0$ for $i<d$ if $H^q_{S_+}(M)_\mu =0$ for $q<d$, proving (i).

For (ii), the condition $H^q_{S_+}(M)_\mu =0$ for $q\leq d$ implies that $H^0_{S_+}(H^d(a;M))_\mu =H^0_{S_+}(H^d(a;M_\mu ))=0$. Hence if $H^d(a;M_\mu )\not= 0$, then $0\not= (S_+)^{\nu -\mu }H^d(a;M)_\mu\subseteq H^d(a;M)_\nu$, which shows (ii) by Lemma \ref{KosCech}. \fini

\bco\label{depthtrunc}
Let $I$ be an $R$-ideal and $M$ be a graded $S$-module. 

Then $\depth_I (M_{\geq \mu})$ is independent of $\mu$ for $\mu>\reg (M)-\depth_{S_+}(M)$.
\eco

\bt\label{assdepth2}
Let $I$ be an $R$-ideal and $M$ be a graded $S$-module with $\reg (M)<\infty$. Set $r:=\reg (M)-\depth_{S_+}(M)$, $d:=\min_{\nu >r}\{ \depth_I (M_\nu )\}=\depth_I (M_{>r})$ and $\mu_0:=\inf\{ \nu >r\ \vert\ \depth_I (M_\nu )=d\}$. 

Then  $\depth_I (M_\mu)=d$ for all $\mu \geq \mu_0$. 
\et

{\it Proof.} We may assume that $d<\infty$.  By definition of $d$, 
$\depth_I (M_\mu )\geq d$ for $\mu \geq \mu_0$, as $\mu_0 >r$. On the other hand, $\depth_I (M_\mu )\leq d$ for 
$\mu \geq \nu$ by Proposition \ref{assdepth} (ii). The conclusion follows.\fini

\section{Cohomological dimension}

Let $A$ be a commutative ring (with unit), $I$ a finitely generated ideal and $M$ a $A$-module. 
First remark that it follows from (\ref{SCL}) that if $M$ is annihilated by an ideal $J$, for instance if 
$J=\ann_A (M)$, considering $M$ as an $A/J$-module one has
\begin{equation}\label{CdAnn}
\cd_I^A (M)=\cd_{(I+J)/J}^{A/J}(M).
\end{equation}
Furthermore,

\bp\label{Cdfgsm}
Let $M$ be a $A$-module.

{\rm (a)} 
 $$
 \cd_I (M)\leq \max_{E\subset M,\ E\ f.g.}\{ \cd_I (E)\} .
 $$
 
{\rm  (b)} 
 $$
 \cd_I (M)\leq   \cd_I (A/\ann_A M)\leq   \cd_I (A).
 $$

{\rm (c)} If $M$ is finitely generated, then
 $$
 \cd_I (M)=\cd_I (A/\ann_A M).
 $$
\ep

{\it Proof.} (a) $M$ is the filtered inductive limit of its submodules of finite type (for the inclusion),
and local cohomology commutes with filtered inductive limits.

(b) The spectral sequence (\ref{SSTor}) shows that $\cd_I (N)\leq \cd_I (A)$ for any module $N$,
in particular $\cd_I (A/\ann_A (M))\leq \cd_I (A)$. Together with (\ref{CdAnn}) applied with $J:=\ann_A (M)$, 
we get $\cd_I^A (M)=\cd_{(I+J)/J}^{A/J}(M)\leq \cd_{(I+J)/J}^{A/J}(A/J)=\cd_I (A/J)$. 

(c) According to (b), it suffices to show that $\cd_I (A/\ann_A M)\leq \cd_I (M)$. Replacing 
$A$ by $A/\ann_A (M)$, we may assume that $M$ is faithful.

We will show that $\cd_I (M)\leq r$ implies $\cd_I (A)\leq r$. This is clear for $r\geq \ara_A (I)$,
and we now perform a descending recursion on $r$. Assume this is true for $r+1$. If
$\cd_I (M)\leq r$, by recursion hypothesis we now that $\cd_I (A)\leq r+1$, hence spectral
sequence (\ref{SSTor}) implies that 
$$
0=H^{r+1}_I (M)=M\otimes_A H^{r+1}_I (A).
$$
As $M$ is faithful and of finite type, \cite[4.3]{St} shows that $H^{r+1}_I (A)=0$,
which implies that $\cd_I (A)\leq r$.
\fini

The following corollary has been proved by M. Dibaei and A. Vahidi in the Noetherian case in
\cite[2.2]{DV}.

\bco\label{dibaei}
Let $M$ be a $A$-module and $I$, $J$ two finitely generated ideals. Then
$$
\cd_{I+J} (M)\leq \cd_I (A/\ann_A M)+\cd_J (M),
$$
and $\cd_{I+J} (M)\leq \cd_I (M)+\cd_J (M)$ if $M$ is finitely generated.
\eco

{\it Proof.} We may assume that $M$ is faithful. If $a$ is a 
finte set of generators of $I$ and $b$ a finite set of generators of $J$, the double 
complex with components $\C^i_a (A)\otimes_A \C^j_b (M)$ gives 
rise to a spectral sequence with second term $H^i_I (H^j_J (M))$ that abuts to 
$H^{i+j}_{I+J}(M)$. Hence,
$$
 \cd_{I+J}(M)\leq \max\{ i+j\ \vert\ H^i_I (H^j_J (M))\not= 0\}\leq \cd_I (A)+\cd_J (M),
$$
by Proposition \ref{Cdfgsm} (b).

Furthermore, by Proposition \ref{Cdfgsm} (c), $\cd_I (A)=\cd_I (M)$ whenever $M$ is faithful and
finitely generated. \fini

\bco\label{cdannihil}
Let $M$ be a finitely generated $A$-module  and $\ia$ an ideal of $A$ such that $\ia^{t}M=0$ 
for some $t$. Then
$$
\cd_I (M)=\cd_I (M/\ia M).
$$
\eco

{\it Proof.} Let $d:=\cd_I (M)=\cd_I (A/\ann_A (M))$ (by Proposition \ref{Cdfgsm}(c)). As $M/\ia M$ is 
annihilated by $\ann_A (M)$, it follows from Proposition \ref{Cdfgsm}(b) that $\cd_I (M/\ia M)\leq d$. 
Furthermore, Proposition \ref{Cdfgsm}(b) shows that the functor $N\mapsto H^d_I (N)$ restricted to
the category of $A$-modules annihilated by $\ann_A (M)$ is right exact. It implies that $H_I^d(M/\ia M)=
H_I^d (M)\otimes_A A/\ia =H_I^d (M)/\ia H_I^d (M)$.

If $\cd_I (M/\ia M)<d$ it implies
$$
H_I^d (M)=\ia H_I^d (M)=\ia^2 H_I^d (M)=\cdots =\ia^t H_I^d (M)=0,
$$
which contradicts the definition of $d$.\fini

\bp\label{cdexseq}
Let  $0\ra M'\ra M\ra M''\ra 0$ be an exact sequence of $A$-modules. Then
$$
\cd_I (M)\leq\max\{ \cd_I (M'),\cd_I (M'')\} \leq \cd_I (A/\ann_A (M)).
$$
Furthermore all inequalities are equalities if $M$ is finitely generated.
\ep

{\it Proof.} First the exact sequences
$$
H^i_I( M')\ra H^i_I( M)\ra H^i_I( M''),\quad i\in {\bf Z},
$$
show the inequality on the left. As both $A$-modules $M'$ and 
$M''$ are annihilated by $\ann_A (M)$, the inequlity on the right follows 
from (\ref{CdAnn}) and Proposition \ref{Cdfgsm}(b).

Finally, the extreme terms are equal according to Proposition \ref{Cdfgsm}(c)
if $M$ is finitely generated.\fini

\brm
If $A$ is a domain, distinct from its field of fractions $K$ and $I$ a proper finitely generated ideal, then $A$ is a 
submodule of $K$ such that 
$$
-\infty =\cd_I (K)<0\leq \cd_I (A).
$$
\erm

\bco\label{cdass}
If $M$ is a Noetherian $A$-module, then
$$
\cd_I (M)=\max_{\ip\in\supp_A (M)}\cd_I (A/\ip )=\max_{\ip\in\ass_A (M)}\cd_I (A/\ip )=\max_{\ip\in\Min_A (M)}\cd_I (A/\ip ).
$$
\eco

{\it Proof.} Let $\Min_A (M)$ be the minimal primes in the support of $M$. Every $\ip\in \supp_A (M)$
contains some $\iq\in \Min_A (M)$, and the canonical epimorphism $A/\iq \ra A/\ip$ gives
an inequality $\cd_I (A/\ip )\leq \cd_I (A/\iq )$ by Proposition \ref{cdexseq}. It follows that
$$
\max_{\ip\in\supp_A (M)}\cd_I (A/\ip )=\max_{\ip\in\Min_A (M)}\cd_I (A/\ip ).
$$
On the other hand $\Min_A (M)\subset \ass_A (M)$ and for $\ip\in \ass_A (M)$ the existence of 
a monomorphism $A/\ip \ra M$ implies by Proposition \ref{cdexseq} that $\cd_I (A/\ip )\leq \cd_I (M)$.
Hence
$$
\max_{\ip\in\Min_A (M)}\cd_I (A/\ip )\leq \max_{\ip\in\ass_A (M)}\cd_I (A/\ip )\leq \cd_I (M),
$$
and it remains to show that $\cd_I (M)\leq \max_{\ip\in\supp_A (M)}\cd_I (A/\ip )$.

Observe that, if $M\not= 0$, it admits a finite filtration by cyclic modules $A/\ip_i$ ($1\leq i\leq t$)
with $\ip_i\in \supp_A (M)$ for all $i$. Again by Proposition \ref{cdexseq}, we obtain
$$
\cd_I (M)=\max_{1\leq i\leq t}\{ \cd_I (A/\ip_i )\}\leq \max_{\ip\in\supp_A (M)}\{ \cd_I (A/\ip )\}.
$$
This concludes the proof.\fini

The following result generalizes the main theorem of \cite{DNT}, which applies in the case of two finitely generated 
modules over a Noetherian ring.

\bp\label{cdsubmod}
Let  $M$ and $N$ be $A$-modules. Assume $M$ is finitely presented and 
$\supp_A (N)\subset \supp_A (M)$. Then
$$
\cd_I (N)\leq \cd_I (M).
$$
\ep

{\it Proof.} Let $E$ be a finitely generated submodule of $N$. The inclusion 
$\supp_A (N)\subset \supp_A (M)$ implies that 
\begin{equation}\label{sqrtann}
\sqrt{\ann_A (M)}\subset \sqrt{\ann_A (E)}.
\end{equation}
As $M$ is finitely generated, the Fitting ideal ${\rm Fitt}^0_A (M)$ has the same radical as $\ann_A (M)$ 
and contains a power of $\ann_A (M)$. Furthermore, as $M$ is finitely presented, this ideal is
finitely generated. It then follows from (\ref{sqrtann}) that there exist $t$ such that
$\ann_A (M)^t\subset \ann_A (E)$. By Corollary \ref{cdannihil} and Proposition \ref{Cdfgsm} (b) and (c), 
one has
$$
\cd_I (E)\leq \cd_I (E/\ann_A (M)E)\leq \cd_I (A/\ann_A (M))=\cd_I (M).
$$
The conclusion follows by Proposition \ref{Cdfgsm} (a) applied to the $A$-module $N$.\fini

Now, let $S$ be a finitely generated standard graded algebra over a commutative ring $R$.

\bd
For a  graded $S$-module $M$, 
$$
a_{S_+}^i (M):=\sup\{ \mu \ \vert\ H^i_{S_+}(M)_\mu \not= 0\}  ,
$$
so that $\reg (M)=\max_i \{ a_{S_+}^i (M)+i\}$.
\ed

\bp\label{stabcd}
Let  $M$ be a finitely generated graded $S$-module and $I$ be a finitely generated $R$-ideal. Then

{\rm (a)} $\cd_I (M_\mu )$ is a non decreasing function of $\mu$ for $\mu > a_{S_+}^0 (M)$,

{\rm (b)} $\cd_I (M_\mu )$ is constant for 
$\mu \geq \reg (M)+n-\depth_\pp (M)$
if $n>0$.
\ep

{\it Proof.} We may, and will, assume that $\cd_\pp (M)>0$, as the proposition is immediate
when $\cd_\pp (M)=0$ by Lemma \ref{3.1}. We also remark that it suffices to prove the claim after the
faithfully flat base change $R\ra R'$, hence we may further assume (making an invertible linear change of coordinates) 
by Remark \ref{glt} that
the sequence $(X_1,\ldots ,X_n)$ is $M$-regular off $V(\pp )$. 

In particular the kernel of the map
$M\ra M(1)$ induced by multiplication by $X_n$ is contained in $H^0_\pp (M)$. 
It follows an
injection $M_\nu \ra M_{\nu+1}$ for $\nu  > a_{S_+}^0 (M)$ which proves (a) by Proposition \ref{cdexseq}
as $M_{\nu+1}$ is finitely generated over $R$.

To prove (b) we consider the two converging spectral sequences arising from the double complex 
$\C^\bullet_I K_\bullet (X_1,\ldots ,X_n ; M)$. They have as respective second terms
$$
\tor_{-p}^{S}(H^q_I (M),R)\quad
{\rm and}
\quad
H^p_I(\tor_{-q}^{S}(M,R)).
$$
Let $d:=\max_{\nu >a_{S_+}^0 (M)}\{ \cd_I (M_\nu )\}$. We may assume $d\geq 0$. 
It follows from the comparison of the spectral sequences that
$$
\tor_{0}^{S}(H^d_I (M),R)_{\nu +1} =H^d_I (M_{\nu +1})/S_1 H^d_I (M_\nu)=0
$$
if $H^{d+i}_I(M_{\nu -i})=0$ for $1\leq i\leq n-1$ and $H^{d+i}_I(\tor_{i}^{S}(M,R)_{\nu +1})=0$ for all $i$. 


But $\tor_i^S (M,R)_{\nu +1}=0$
for $\nu \geq \reg (M)+i$ and $\tor_{i}^{S}(M,R)=0$ for $i>n-\depth_\pp (M)$ by Lemma \ref{KosCech}.

It follows that $H^d_I (M_\nu)=0$ implies $H^d_I (M_{\nu +1})=0$ if 
 $$
 \nu \geq \max\{ a_{S_+}^0 (M)+n,\reg (M)+n-\depth_\pp (M)\} =\reg (M)+n-\depth_\pp (M) .
 $$
This implies (b), in view of (a).
 \fini

\section{Associated primes of the graded components of a graded module}

Let $S$ be a standard graded Noetherian algebra over a commutative ring $R$.

\bt\label{assgraded}
Let $M$ be a graded $S$-module. 
then
 $$
 \bigcup_{\mu \in \ZZ}\ass_R (M_\mu )=\{ \Ip \cap R,\ \Ip\in \ass_S (M)\} .
 $$
\et

{\it Proof.} 
For $\mu\in\ZZ$, let $\ip\in \ass_R (M_\mu )$. There exists $x\in M_\mu$ with 
$\ip =\ann_R (x)$. Hence $\ip R_\ip =\ann_{R_\ip}(x)$. Let $\Iq$ be a $S \otimes_R R_\ip$-ideal,
maximal among those of the form $\ann_{S \otimes_R R_\ip}(y)$, $y\in M\otimes_R R_\ip$,
that contains    $\ann_{R_\ip}(x)$. The ideal $\Iq$ is associated to $M\otimes_R R_\ip$, hence 
$\Ip:=\Iq \cap S$ is associated to $M$ and $\Ip \cap R=\ip$. One inclusion follows.

Conversely, let $\Ip$ be an ideal associated to $M$. We need to show that $\ip :=\Ip\cap R$ is associated 
to $M_\mu$ for some $\mu$. This will be the case if $\ip R_\ip$ is associated $(M_\mu )_\ip$, so that we may assume that $R$ is local with maximal ideal $\ip$. Let $m\not= 0$ in $M$ such that $\Ip m=0$. If 
$m_\nu$ is the degree $\nu$ component of $m$, one has $\ip m_\nu=0$. Hence choosing $\mu$ such
that $m_\mu \not= 0$, one has $\ip \subseteq \ann_R (m_\mu )$, hence $\ip =\ann_R (m_\mu )$, as $\ip$ is maximal.\fini

\bt\label{assymass}
Let $M$ be a graded $S$-module. If $H^0_{S_+}(M)_\nu =0$, then $\ass_R (M_\nu )\subseteq \ass_R (M_{\mu})$ for all $\mu\geq \nu$.
\et

{\it Proof.} Let $x_1,\ldots ,x_n$ be generators of $S_1$ as an $R$-module and $T_1,\ldots ,T_n$ be variables. Set $P:=\sum_{|\a |=\mu -\nu}
x^{\a}T^{\a}$. The polynomial $P$ is of bidegree $(\mu -\nu ,\mu -\nu )$ for the bigraduation defined by setting $\deg (x_i):=(1,0)$ and $\deg (T_i):=(0,1)$, and $c(P)=(S_+)^{\mu -\nu}$.  Theorem \ref{LDM} shows that the kernel $K$ of the map 
$$
\xymatrix{M[\ult ]\ar^{\times P}[r]&M[\ult ]\\}
$$
is a submodule of $H^0_{S_+}(M)[\ult ]$. Hence $K$ vanishes in bidegree $(\nu ,\theta )$ for any $\theta$. In particular, it provides an injective map
$$
 \xymatrix{M_\nu =M[\ult ]_{\nu ,0}\ar^{\times P}[r]&M[\ult ]_{\mu ,\mu -\nu}.\\}
$$
As $M[\ult ]_{\mu ,\mu -\nu}$ is a finite direct sum of copies of $M_\mu$, it follows that any associated prime of $M_\nu$
is associated to $M_\mu$.
\fini

\bt\label{assympolass}
Let
$M$ be a finitely generated graded $S$-module and $\A$ be the finite set $\cup_{\mu\in \ZZ}\ass_R (M_\mu )$ (Theorem \ref{assgraded}). 
Set 
$$
j(M):=\max_{\ip\in \A}\{ a_{S_+}^{0} (H^0_\ip (M\otimes_{R}R_\ip))\}\leq a_{S_+}^{0}(M).
$$

 Then
 for any ideal $\ip\in\spec (R)$, 
 $\ell_{R_\ip}(H^0_\ip (M_\mu\otimes_{R}R_\ip))$ is a nondecreasing function of $\mu$, for $\mu>j(M)$. \et

{\it Proof.} First notice that $H^0_\ip (M\otimes_{R}R_\ip)=0$ if $\ip\not\in \A$. 
Let $\ip \in \A$. We will prove that $\ell_{R_\ip}(H^0_\ip (M_\mu \otimes_{R}R_\ip))$ is a nondecreasing function of $\mu$, for $\mu>a_{S_+}^{0} (H^0_\ip (M\otimes_{R}R_\ip))$. The proof of Theorem \ref{assymass},  applied to the $S\otimes_R R_\ip$-module $H^0_\ip (M\otimes_{R}R_\ip)$, with 
$P:=\sum x_iT_i$ provides an injective morphism of $R_\ip [\ult ]$-modules
$$
 \xymatrix{H^0_\ip (M\otimes_{R}R_\ip)_{\mu}[\ult ]\ar[r]& H^0_\ip (M\otimes_{R}R_\ip)_{\mu +1}[\ult ] .\\}
$$

Let $R_\ip (\ult ):=S^{-1}R_\ip [\ult ]$ with $S$ the multiplicative system of polynomials whose coefficient ideal  is the unit ideal. The above injection induces
an injective morphism of $R_\ip (\ult )$-modules of finite length
$$
 \xymatrix{H^0_\ip (M\otimes_{R}R_\ip)_{\mu}\otimes_{R_\ip }R_\ip(\ult )\ar[r]& H^0_\ip (M\otimes_{R}R_\ip)_{\mu +1}\otimes_{R_\ip }R_\ip(\ult ) .\\}
$$
For any $R_\ip$-module $N$ of finite length, the $R_\ip(\ult )$-module $N\otimes_{R_\ip }R_\ip(\ult )$ is a module of the same length as the
$R_\ip$-module $N$. The conclusion follows.
\fini

\bco\label{tametop}
Let 
$M$ be a finitely generated graded $S$-module. Then for any $\ip\in\spec (R)$ of height 0, 
 $\ell_{R_\ip}(M_\mu\otimes_{R}R_\ip )$ is a nondecreasing function of $\mu$, for $\mu>a^0_{S_+} (M)$.
 In particular $M_\mu \otimes_R R_\ip =0$ for all $\mu>\reg (M)$ or $M_\mu \otimes_R R_\ip \not= 0$ for all $\mu>\reg (M)$.
 \eco
 
{\it Proof.} Notice that $M\otimes_{R}R_\ip =H^0_\ip (M\otimes_{R}R_\ip)$ as $\ip$ is of height 0.
Also recall that $M$ (hence $M\otimes_{R}R_\ip$) is generated in degrees at most $\reg (M)$.\fini

\bl\label{qsM}
Let $M$ be a graded $S$-module generated in degree at most $B$. Then,

(i)  $M_{\mu} \not= 0\ \Rightarrow M_{\mu +1}\not= 0$, if  $\mu >a_{S_+}^0 (M)$,

(ii)  $M_{\mu} \not= 0\ \Leftrightarrow M_{\mu +1}\not= 0$ if  $\mu >\max\{ B,a_{S_+}^0 (M)\}$, 

(iii)  $M_{\mu} \not= 0\ \Leftrightarrow M_{\mu +1}\not= 0$ if  $\mu >\reg (M)$.
\el

{\it Proof.} $(i)$ follows from Theorem \ref{assymass} as $R$ is Noetherian.  Now $(ii)$ and $(iii)$ follows from $(i)$,
as $M_\mu = 0\Rightarrow M_{\mu +1}= 0$ for $\mu\geq B$ and $M$ is generated in degrees at most $\reg (M)$. \fini

\brm
The proof of Theorem \ref{assymass} shows that the above lemma holds without assuming that $R$
is Noetherian.
\erm

\section{Duality results}

\subsection{Preliminaries on RHom}

Let $(X,\OO_X )$ be a ringed space and $Y$ be closed in $X$. For any complex $K^\bullet \in D(X)$, let 
$C(K^\bullet )$ be the Godement resolution of $K^\bullet$, and set 

-- $\RGG (X,K^\bullet ):=C(K^\bullet )$ and $\RGG_Y (X,K^\bullet ):= \GG_{Y} (X,C(K^\bullet ))$ in $D(X)$,

-- $\RG (X,K^\bullet ):=\G (X,C(K^\bullet ))$ and $\RG_Y (X,K^\bullet ):= \G_{Y} (X,C(K^\bullet ))$ in $D(\G (X,\OO_X ))$.

Notice that a flasque resolution of $K^\bullet$  in $D^+ (X)$ (e.g. an injective resolution) can be used in place of $C(K^\bullet )$
if $K^\bullet\in D^+ (X)$.

We set $H^i (X,K^\bullet ):=H^i (\RG (X,K^\bullet ))$, $H^i_Y (X,K^\bullet ):=H^i (\RG_Y (X,K^\bullet ))$ and ${\bf H}^i (X,K^\bullet ):=H^i (\RGG (X,K^\bullet ))$, which coincides with the usual notations for $\OO_X$-modules when considered as complexes concentrated in degree 0.

If there exists $d$ such that, for any $\OO_X$-module $\E$, $H^i (X,\E )=0$ (resp. $H^i_Y (X,\E )=0$) for $i>d$ then any flasque resolution of
$K^\bullet$ can be used in place of $C(K^\bullet )$ to compute $H^i (X,K^\bullet )$ and ${\bf H}^i (X,K^\bullet )$ (resp. $H^i_Y (X,K^\bullet )$).

Given  $K^\bullet$  in $D(X)$ and  $L^\bullet$  in $D^+ (X)$ one checks that the class in $D(X)$ (resp. in $D(\G (X,\OO_X ))$)
of  the complex $\shhom_{\OO_X}^\bullet (K^\bullet ,I^\bullet )$ (resp. $\hom_{\OO_X}^\bullet (K^\bullet ,I^\bullet )$) is independant of the 
choice of an injective resolution $I^\bullet$ of $L^\bullet$ in $D^+ (X)$, and one set
$$
R\!\hom_{\OO_X}^\bullet (K^\bullet ,L^\bullet ):=\hom_{\OO_X}^\bullet (K^\bullet ,I^\bullet )\ {\rm and}\ 
R\!\shhom_{\OO_X}^\bullet (K^\bullet ,L^\bullet ):=\shhom_{\OO_X}^\bullet (K^\bullet ,I^\bullet ).
$$
When the components of the complex $K^\bullet$ are locally free $\OO_X$-modules,
there is a quasi-isomorphism
$$
R\!\shhom_{\OO_X}^\bullet (K^\bullet ,L^\bullet )\simeq\shhom_{\OO_X}^\bullet (K^\bullet ,L^\bullet ).
$$
For a pair $(\E ,\F )$ of $\OO_X$-modules, one defines
$$
\hom_Y (\E ,\F ):=\G_Y (X,\hom_{\OO_X}(\E ,\F ))=\hom_{\OO_X}(\E ,\GG_Y \F )
$$
and
$$
\shhom_Y (\E ,\F ):=\GG_Y \shhom_{\OO_X}(\E ,\F )=\shhom_{\OO_X}(\E ,\GG_Y \F ),
$$
and then extends these definitions to pairs of complexes of $\OO_X$-modules
as usual.

Assume $L^\bullet$ is bounded below. Given a bounded below injective resolution $I^\bullet$
of $L^\bullet$, the components of the complex $\shhom_{\OO_X}^\bullet (K^\bullet ,I^\bullet )$ are 
flasques hence 
$$
\begin{array}{rl}
\RGG_Y (R\!\shhom_{\OO_X}^\bullet (K^\bullet ,L^\bullet ))&=\RGG_Y (\shhom_{\OO_X}^\bullet (K^\bullet ,I^\bullet ))\\
&=\GG_Y (\shhom_{\OO_X}^\bullet (K^\bullet ,I^\bullet ))\\
&=\shhom_{Y}^\bullet (K^\bullet ,I^\bullet ),\\
\end{array}
$$
in the following cases :

(a) $K^\bullet$ is bounded below,

(b) there exists $d$ such that, for any $\OO_X$-module $\E$, $H^i_Y (X,\E )=0$ for $i>d$.

In cases (a) and (b) $\shhom_{Y}^\bullet (K^\bullet ,I^\bullet )$ is independent of the choice of $I^\bullet$ (up
to an isomorphism in $D(X)$), and one sets 
$$
\RGG_Y^\bullet (K^\bullet ,L^\bullet ):=\shhom_{Y}^\bullet (K^\bullet ,I^\bullet ).
$$

As the components of $\GG_Y I^\bullet$ are injective $\OO_X$-modules, 
$\shhom_{Y}^\bullet (K^\bullet ,I^\bullet )\simeq \shhom_{\OO_X}^\bullet (K^\bullet ,\GG_Y I^\bullet)
\simeq R\!\shhom_{\OO_X}^\bullet (K^\bullet ,\RGG_Y L^\bullet )$. In other words,
$$
\RGG_Y^\bullet (K^\bullet ,L^\bullet )\simeq R\!\shhom_{\OO_X}^\bullet (K^\bullet ,\RGG_Y L^\bullet )
$$
if (a) or (b) holds.

\subsection{Some spectral sequences}

Let $A$ be a commutative ring (with unit) and $I$ be a finitely generated $A$-ideal. Set 
$(X,\OO_X ):=(\spec (A),\tilde A )$ and $Y:=V(I)\subset X$. 

Assume that there exists $n$ such that 
$$
H^i_I (A)=H^i_ Y (X,\OO_X )=0,\ {\rm for}\ i\not= n.
$$
Then $\RGG_Y (\OO_X )\simeq {\bf H}^n_Y (\OO_X )[-n]$ in $D^+ (X,\OO_X )$. Given a complex $K^\bullet \in 
D^- (X,\OO_X )$, it follows that 
$$
\begin{array}{rl}
\RGG_Y (R\!\shhom_{\OO_X}^\bullet (K^\bullet ,\OO_X ))&\simeq R\!\shhom_{\OO_X}^\bullet (K^\bullet ,\RGG_Y (\OO_X ))\\
&\simeq R\!\shhom_{\OO_X}^\bullet (K^\bullet ,{\bf H}^n_Y (\OO_X ))[-n]\\
\end{array}
$$
in $D(X,\OO_X )$. Such an isomorphism holds for $K^\bullet \in D(X,\OO_X )$ when $X$ has finite homological dimension, hence, for instance, if $A$ is Noetherian of finite dimension.

Assuming further that the components of $K^\bullet$ are locally free $\OO_X$-modules of finite type, under one of the two hypotheses above, one has
$$
\RGG_Y (\shhom_{\OO_X}^\bullet (K^\bullet ,\OO_X ))\simeq \shhom_{\OO_X}^\bullet (K^\bullet ,{\bf H}^n_Y (\OO_X ))[-n]\ {\rm in}\ D(X,\OO_X ).
$$
This provides a spectral sequence
$$
E_2^{p,q}={\bf H}^p_Y (H^q \shhom_{\OO_X}^\bullet (K^\bullet ,\OO_X ))\Rightarrow H^{p+q-n}
( \shhom_{\OO_X}^\bullet (K^\bullet ,{\bf H}^n_Y (\OO_X ))),
$$
in the two cases above.

As the $\OO_X$-modules taking place in this spectral sequence are quasi-coherent, and the components of $K^\bullet$ are 
finitely presented (recall that $\widetilde{\hom_A (M,N)}=\shhom_{\OO_X}(\tilde{M},\tilde{N})$ if $M$ is finitely presented) the above 
results shows the following Proposition.

\bp
Let $A$ be a commutative ring and $I$ a finitely generated $A$-ideal. Assume that there exists $n$ such that
$$
H^i_I (A)=H^i_Y (X,\OO_X )=0\ {\rm for}\ i\not= n.
$$
Then, for any bounded below complex of projective $A$-modules of finite type (resp. for any  complex of projective $A$-modules of finite type if $A$ is Noetherian of finite dimension) $K^\bullet$, there exists a spectral sequence
$$
E_2^{p,q}=H^p_I (H^q \hom_{A}^\bullet (K^\bullet ,A))\Rightarrow H^{p+q-n}
( \hom_{A}^\bullet (K^\bullet ,H^n_I (A))).
$$
\ep

\bco
Assume $(A,\im )$ is local Gorenstein of dimension $n$. Then for any complex $K^\bullet$ of free $A$-modules of
finite type, there is a spectral sequence
$$
E_2^{p,q}=H^p_\im (H^q \hom_{A}^\bullet (K^\bullet ,A))\Rightarrow H^{p+q-n}
( \hom_{A}^\bullet (K^\bullet ,H^n_\im (A))).
$$
\eco

{\it Proof.} Under the hypotheses of the Corollary, $H^i_\im (A)=0$ for $i\not= n$ and
$H^n_\im (A)$ is injective.\fini

\bex
In the context of the Corollary, taking for $K^\bullet$ the dual of a resolution of finitely generated module $M$ by free modules 
of finite type it gives the local duality 
$$
H^p_\im (M)\simeq \hom_A (\ext^{n-p}_A (M,A),H^n_\im (A)).
$$
\eex

\bco\label{SSwregseq}
Assume $I$ is generated by a weakly regular sequence of length $n$, then for any bounded below complex $K^\bullet$ of projective$A$-modules of
finite type, there is a spectral sequence
$$
E_2^{p,q}=H^p_I (H^q \hom_{A}^\bullet (K^\bullet ,A))\Rightarrow H^{p+q-n}
( \hom_{A}^\bullet (K^\bullet ,H^n_I (A))).
$$
\eco

\bex
Let $R$ be a commutative ring, $X_i$ for $1\leq i\leq n$ be indeterminates, set $A:=R[X_1,\ldots ,X_n]$, with
its standard grading and $\ip :=(X_1,\ldots ,X_n)$. As $H^n_\ip (A)\simeq (X_1\cdots X_n)^{-1}R[X_1^{-1},\ldots ,X_n^{-1}]$, it follows that for any graded free $A$-module of finite type $F$ and every integer $\nu$, the pairing
$$
\xymatrix@R=3pt{
\hom_A (F,H^n_\ip (A))_{-\nu -n}\otimes_R F_\nu\ar[r]&H^n_\ip (A)_{-n}\simeq R\\
(u:F\ra H^n_\ip (A)(-\nu -n))\otimes_R x\ar@{|-{>}}[r]&u(x)\\}
$$
defines a perfect duality between $R$-modules of finite type, and this duality is functorial in the free graded $A$-module $F$. It gives, for each integer $\nu$, an isomorphism of complexes of $R$-modules 
$$
\hom_A^\bullet (K^\bullet ,H^n_\ip (A))_{-\nu -n}\simeq \hom_R^\bullet (K^\bullet_\nu,R).
$$
Together with Corollary \ref{SSwregseq}, it gives for any $\nu$ a spectral sequence of $R$-modules
$$
E_2^{p,q}=H^p_\ip (H^q \hom_{A}^\bullet (K^\bullet ,A))_\nu \Rightarrow H^{p+q-n}
( \hom_{R}^\bullet (K^\bullet_{-\nu -n},R)).
$$

Replacing $K^\bullet$ by its dual $F^\bullet$, one deduces of a spectral sequence
$$
E_2^{p,q}=H^p_\ip (H^q (F^\bullet )_\nu ) \Rightarrow H^{p+q-n}
( \hom_{R}^\bullet (\hom_{A}^\bullet (F^\bullet ,A)_{-\nu -n},R)).
$$
For instance, if $M$ is a graded $A$-module admitting a resolution by free modules of finite rank, taking for $K^\bullet$
such a resolution, that one may assume to be graded, it follows a spectral sequence :
$$
E_2^{p,q}=H^p_\ip (\ext^q_A (M,A))_\nu  \Rightarrow \ext^{p+q-n}_R
(M_{-\nu -n},R).
$$

\eex

\subsection{The Herzog-Rahimi spectral sequences}

We keep notations as in the preceding subsection. For any graded complex $K^\bullet$ whose components are
of finite type, we have established an isomorphism
$$
\RGG_\ip (\hom_{A}^\bullet (K^\bullet ,A ))\simeq \hom_{A}^\bullet (K^\bullet ,H^n_\ip (A))[-n]\ {\rm in}\ D(A),
$$
whenever $K^\bullet$ is bounded above or $A$ is Noetherian of finite dimension. 

In each of these cases, it follows, for any integer $\nu$, isomorphisms
$$
\RGG_\ip (\hom_{A}^\bullet (K^\bullet ,A ))_\nu \simeq \hom_{R}^\bullet (K^\bullet_{-\nu -n} ,R)[-n]\ {\rm in}\ D(R).
$$
Hence if $F^\bullet$ is a graded complex of finite free $A$-modules, and either $A$ is Noetherian of finite
dimension or $F^\bullet$ is bounded below, one has isomorphisms
$$
\begin{array}{rl}
\RGG_\ip (F^\bullet)_\nu &\simeq \hom_{R}^\bullet (\hom_{A}^\bullet (F^\bullet ,A )_{-\nu -n} ,R)[-n]\\
&\simeq R\!\hom_{R}^\bullet (\hom_{A}^\bullet (F^\bullet ,A )_{-\nu -n} ,R)[-n]\ {\rm in}\ D(R).\\
\end{array}
$$
Now, assume further that $(R,\im )$ is local Gorenstein of dimension $d$. The above isomorphisms then give
$$
\begin{array}{rl}
R\!\hom_{R}^\bullet (\RGG_\ip (F^\bullet)_\nu ,H^d_\im (R))
&\simeq  R\!\hom_{R}^\bullet (\RGG_\ip (F^\bullet)_\nu ,\RGG_\im (R))[d]\\
&\simeq  \RGG_\im (R\!\hom_{R}^\bullet (\RGG_\ip (F^\bullet)_\nu , R))[d]\\
& \simeq \RGG_\im (R\!\hom_{R}^\bullet (R\!\hom_{R}^\bullet (\hom_{A}^\bullet (F^\bullet ,A )_{-\nu -n} ,R) ,R))[n+d]\\
& \simeq \RGG_\im (\hom_{A}^\bullet (F^\bullet ,A )_{-\nu -n})[n+d]\\
\end{array}
$$

As $R$ is Gorenstein, $H^d_\im (R)$ is the injective envelope of the residue field of $R$, and we obtain a
spectral sequence
$$
E_2^{p,q}=H^p_\im (H^q (\hom_{A}^\bullet (F^\bullet ,A )_{-\nu -n}) ) \Rightarrow
\hom_{R} ( H^{n+d-p-q}(\RGG_\ip (F^\bullet)_\nu ),H^d_\im (R)).
$$
If $M$ is graded $A$-module with free resolution $F^\bullet$, this spectral sequence takes the 
form
$$
E_2^{p,q}=H^p_\im (\ext^q_A (M ,A )_{-\nu -n}) \Rightarrow
\hom_{R}( H^{n+d-p-q}_\ip (M)_\nu ,H^d_\im (R)).
$$
Which is the Herzog-Rahimi spectral sequence, as $\om_A \simeq A(-n)$ in this situation.

\section{Tameness of local cohomology over Noetherian rings}

In this section $S$ is a finitely generated standard graded algebra over an epimorphic image $R$ of
a Gorenstein ring.

\bt\label{duality}
Let $(R,\im )$ be a local Noetherian Gorenstein ring of dimension $d$, $S$ be a finitely generated standard graded 
Cohen-Macaulay algebra over $R$ and $M$ be a finitely
generated graded $S$-module. Set $\hbox{---}^{\vee}:=\hom_R (\hbox{---},H^d_\im (R))$. Then, there is a spectral sequence
$$
E^{i,j}_{2}=H^{i}_{\im}(\Ext^{j}_{S}(M,\omega_{S})_{\g})\Rightarrow
(H^{\dim S-(i+j)}_{S_{+}}(M)_{-\g})^{\vee}.
$$
\et

{\it Proof.} See [R, section 4] or the previous section.\fini\medskip

For a $R$-module or a $S$-module $M$, we set 
$$
H^0_{[i]}(M):=\bigcup_{I\subseteq R,\; \dim (R/I)\leq i}H^0_{I}(M)\subseteq M.
$$

We will need the following facts about the functor $H^0_{[i]}(\hbox{---})$,

\bl\label{localcohsuppdim}
Let $M$ be a graded $S$-module. Then

(i) $H^0_{[i]}(M)_\g =H^0_{[i]} (M_\g )$.

(ii) If $\ip$ a prime ideal of  $R$ with $\dim (R/\ip)=i$, then
$$
H^0_{[i]}(M)\otimes_{R}R_\ip =H^0_\ip (M\otimes_{R}R_\ip).
$$
\el

{\it Proof.} Claim (i) follows from the fact that, for any $R$-ideal $I$, $H^0_I(M)_\g =H^0_I (M_\g)$.
For (ii), recall that inductive limits commute with tensor products and notice that if $\dim (R/I)\leq i$,
 then $H^0_{I_\ip} (M\otimes_{R}R_\ip)=0$ if $I\not\subseteq \ip$, and 
$H^0_{I_\ip} (M\otimes_{R}R_\ip)=H^0_\ip (M\otimes_{R}R_\ip)$ else.\fini

\bt\label{tame}
Let $S$ be a polynomial ring in $n$ variables over an equidimensional  Gorenstein ring $R$ of dimension $d$. 
Let $M$ be a finitely generated graded $S$-module.

Then there exists $A,B,C,D,E$ defined below such that :

(1a)  For  $\g >A$,
$
\dim H^i_{S_+}(M)_{-\g}=d \Rightarrow \  \dim H^i_{S_+}(M)_{-\g-1}=d.
$

(1b)  For  $\g >B$,
$
\dim H^i_{S_+}(M)_{-\g}=d \Leftrightarrow \  \dim H^i_{S_+}(M)_{-\g-1}=d.
$

(2a)  For  $\g >C$,
$
\dim H^i_{S_+}(M)_{-\g}\geq d-1 \Rightarrow \  \dim H^i_{S_+}(M)_{-\g-1}\geq d-1.
$

(2b)  For  $\g >D$,
$
\dim H^i_{S_+}(M)_{-\g}\geq d-1 \Leftrightarrow \  \dim H^i_{S_+}(M)_{-\g-1}\geq d-1.
$

(3) For  $\g>E$,
$
\dim H^i_{S_+}(M)_{-\g}\geq d-2 \Rightarrow \  \dim H^i_{S_+}(M)_{-\g-1}\geq d-2.
$

In particular, there exists $\g_0$ such that either $\dim H^i_{S_+}(M)_{\g}$
is constant for $\g <\g_0$ of value at least $d-2$, or $\dim H^i_{S_+}(M)_{\g}<d-2$
for any $\g <\g_0$.

Set $a^i_j :=\fin (H_{S_+}^{0} (H^0_{[d-j]}(\ext^{i}_S(M,\om_S ))))\leq \fin (H^0_{S_+}(\ext^{i}_S(M,\om_S )))$ and $r^{i}_{j}:=\reg (H^0_{[d-j]}(\ext^{i}_S(M,\om_S )))$. Then one has :

$A:=a_{0}^{n-i}$, $B:=r^{n-i}_{0}$, $C:=\max\{a_{1}^{n-i+1},a_{0}^{n-i}\}$, 
$D:=\max\{r_{1}^{n-i+1},r_{0}^{n-i}\}$ and $E:=\max\{ a_{0}^{n-i},r_{0}^{n-i+1}-1,r_{2}^{n-i+1}-2,a_{2}^{n-i+2}\}$.

\et
{\it Proof.} Recall that if $N$ is a finitely generated $R$-module, $\dim N<r$ if and only if
$N_\ip =0$ for all $\ip\in \spec (R)$ with $\dim (R/\ip )=r$. Furthermore, it follows from
Lemma \ref{localcohsuppdim}  and from
Lemma \ref{localreg} and its proof that for any $\ip\in \spec (R)$ with $\dim (R/\ip )=r$ and for any $\ell$, one has
$$
a^\ell_{S_+}(H^0_{\ip}(\ext^{i}_{S\otimes_{R}R_\ip}(M\otimes_{R}R_\ip,\om_{S\otimes_{R}R_\ip} )))
\leq a^\ell_{S_+}(H^0_{[r]}(\ext^{i}_S(M,\om_S ))).
$$
Therefore, 
$
\reg (H^0_{\ip}(\ext^{i}_{S\otimes_{R}R_\ip}(M\otimes_{R}R_\ip,\om_{S\otimes_{R}R_\ip} )))\leq r^i_{n-r}
$,
and, as a consequence, it suffices to prove $(1a)$ and $(1b)$
when $\dim R=0$, $(2a)$ and $(2b)$ when $\dim R=1$, and $(3)$ when $\dim R=2$.

Notice that $H^0_{[d]}(N)=N$ for any $S$-module $N$.

 If $\dim R=0$, by Theorem \ref{duality}, 
$(H^{i}_{S_{+}}(M)_{-\g })^\vee\simeq \Ext^{n-i}_{S}(M,\omega_{S})_{\g }$, and the result
follows from Lemma \ref{qsM}(i) and (iii).

If $\dim R=1$,  Theorem \ref{duality} provides exact sequences
$$
0\ra H^1_\im ( \Ext^{n-i}_{S}(M,\omega_{S})_\g )
 \ra (H^{i}_{S_+}(M)_{-\g})^\vee\ra H^0_{\im}( \Ext^{n-i+1}_{S}(M,\omega_{S})_\g )
\ra 0.
$$

The result follows applying Lemma \ref{qsM} (i) and (iii) to $H^0_{\im}( \Ext^{n-i+1}_{S}(M,\omega_{S}))$
 and  Corollary \ref{tametop} to $\Ext^{n-i}_{S}(M,\omega_{S})_\g$. 
 Indeed, $H^1_\im ( \Ext^{n-i}_{S}(M,\omega_{S})_\g)$
is zero if and only if $\dim \Ext^{n-i}_{S}(M,\omega_{S})_\g <1$. 

We now assume that $\dim R =2$. In this case Theorem \ref{duality} provides a spectral sequence 
which converges to $(H^{\bullet}_{S_+}(M)_{-\g})^\vee$, 

$$
\small{
\xymatrix{
\cdots&H^0_{\im}(\ext^{n-i+1}_{S}(M,\om_S)_{\g})\ar_{\psi^{n-i+1}_{\g}}[ddl]&
H^0_{\im}(\ext^{n-i+2}_{S}(M,\om_S)_{\g})\ar_(.4){\psi^{n-i+2}_{\g}}[ddl]\\
\cdots&H^1_{\im}(\ext^{n-i+1}_{S}(M,\om_S)_{\g})&H^1_{\im}(\ext^{n-i+2}_{S}(M,\om_S)_{\g})\\
H^2_{\im}(\ext^{n-i}_{S}(M,\om_S)_{\g})&H^2_{\im}(\ext^{n-i+1}_{S}(M,\om_S)_{\g})&\cdots\\
}}
$$
It provides a filtration $F_{*}^0\subseteq F_{*}^1\subseteq F_{*}^2 =(H^{i}_{S_+}(M)_{-*})^\vee$,
by graded $S$-modules,
such that $F_\g^2/F_\g^1\simeq \ker (\psi^{n-i+2}_{\g})$, $F_\g^1/F_\g^0\simeq H^1_{\im}(\ext^{n-i+1}_{S}(M,\om_S)_{\g})$ and $F_\g^0\simeq \coker (\psi^{n-i+1}_{\g})$.

We will show that,  the three modules satisfy :

(i)  $F^0_\g \not= 0\Rightarrow F^0_{\g +1} \not= 0$, if $\g >a_{0}^{n-i}$,
 
(ii)  $F^1_\g /F^0_\g \not= 0\Rightarrow F^1_{\g +1} /F^0_{\g +1} \not= 0$ if $\g >\max\{ r_{0}^{n-i+1}-1,r_{2}^{n-i+1}-2\}$,
 
(iii) $F^2_\g /F^1_\g \not= 0\Rightarrow F^2_{\g +1} /F^1_{\g +1} \not= 0$ if $\g >a_{2}^{n-i+2}$.

For (i) notice that $\coker (\psi^{n-i+1}_{\g})=0$ if and only if $H^2_{\im}(\ext^{n-i}_{S}(M,\om_S)_{\g})=0$,
hence  if and only if $\dim (\ext^{n-i}_{S}(M,\om_S)_{\g})<2$. Hence (i) follows from Corollary \ref{tametop}.

For (ii), let $N:=\ext^{n-i+1}_{S}(M,\om_S)/
 H^0_{\im}(\ext^{n-i+1}_{S}(M,\om_S))$. The exact sequence
 $$
 0\ra  H^0_{\im}(\ext^{n-i+1}_{S}(M,\om_S))\ra \ext^{n-i+1}_{S}(M,\om_S)\ra N\ra 0
$$
 shows that $H^0_\im (N)=0$ (hence $\depth (N)\geq 1$),
 $$
 F^1/F^0=H^1_{\im}(\ext^{n-i+1}_{S}(M,\om_S))\simeq H^1_{\im}(N)
 $$
and 
$$
\begin{array}{rl}
a_{S_+}^{1}(N)&\leq \max \{ a_{S_+}^{1}(\ext^{n-i+1}_{S}(M,\om_S)),a_{S_+}^{2}(H^0_{\im}(\ext^{n-i+1}_{S}(M,\om_S)))\}\\
&\leq  \max\{ r_{0}^{n-i+1}-1,r_{2}^{n-i+1}-2\} .\\
\end{array}
$$
Hence, Proposition \ref{assdepth} (ii) implies (ii)\smallskip

For (iii), let  $\a \in H^0_\im (\ext^{n-i+2}_{S}(M,\om_S)_{\g})= H^0_\im (\ext^{n-i+2}_{S}(M,\om_S))_{\g}$. 

Let $x_1,\ldots ,x_t$ be generators of $S_1$ as an $R$-module and set  $\ell :=\sum_i x_iT_i\in S[\ult ]$. For $\g>a_{S_+}^0 (H^0_\im (\ext^{n-i+2}_{S}(M,\om_S)))$,  $\ell\a$ is not zero in $H^0_\im (\ext^{n-i+2}_{S}(M,\om_S)_{\g +1})[\ult ]$ by Theorem \ref{LDM} (see also the proof of Theorem \ref{assympolass}). The commutative diagram
$$
\xymatrix{
H^0_{\im}(\ext^{n-i+2}_{S}(M,\om_S)_{\g})
\ar^(.45){\times \ell}[r]\ar_{\psi^{n-i+2}_{\g}}[d]
&H^0_{\im}(\ext^{n-i+2}_{S}(M,\om_S)_{\g +1})[\ult ]
\ar^{\psi^{n-i+2}_{\g +1}\otimes_S 1_{S[\ult ]}}[d]\\
H^2_{\im}(\ext^{n-i+1}_{S}(M,\om_S)_{\g})
\ar_(.45){\times \ell}[r]&
H^2_{\im}(\ext^{n-i+1}_{S}(M,\om_S)_{\g +1})[\ult ]\\
}
$$
then shows that $\psi^{n-i+2}_{\g}(\a )\not= 0$ if $\psi^{n-i+2}_{\g +1}$ is injective. 
Hence, $\psi^{n-i+2}_{\g}$ is injective if $\psi^{n-i+2}_{\g +1}$ is. Claim (iii) follows.
\fini
\bigskip\bigskip

\bt\label{tamedim2}
Let $S$ be a Noetherian standard graded algebra over a commutative ring $R$.
Assume $R$  has dimension at most two and either $R$ is an epimorphic image of a Gorenstein ring or 
$R$ is local.
Let $M$ be a finitely generated graded $S$-module.

Then there exists $\g_0$ such that, for any $i$, 

\centerline{\{$H^i_{S_+}(M)_\g =0$ for $\g <\g_0$\} or 
\{$H^i_{S_+}(M)_\g \not= 0$ for $\g <\g_0$\}.}
\et

{\it Proof.} First, if $R$ is local, then we can complete $R$ to reduce to the case where
$R$ is a quotient of a regular ring (by Cohen structure theorem), hence an epimorphic image of a 
Gorenstein ring.

As a Gorenstein ring is a finite product of equidimensional Gorenstein rings,
and each such ring is itself a quotient of a Gorenstein ring of any bigger dimension, $R$ is
also a quotient of an equidimensionnal Gorenstein ring $R'$. We further remark that $R$ 
is the epimorphic image of $R'/K$, where $K$ is generated by a regular sequence of length 
$\dim R'-2$ in $R'$. 

Thus we may, and will, assume that $R$ is an equidimensional Gorenstein ring of dimension 
at most two. Now $S$ is an epimorphic image of a polynomial ring in a finite number of variables 
over $R$, so that we may, and will, also assume that $S$ is a polynomial ring over $R$.

The result then follows from Theorem \ref{tame}.\fini


\begin{thebibliography}{1}

\bibitem{Br1}
Markus Brodmann.
\newblock Asymptotic stability of Ass($M/I^nM$). 
\newblock Proc. Amer. Math. Soc. {\bf 74} (1979), no. 1, 16-18.

\bibitem{Br2}
Markus Brodmann.
\newblock A cohomological stability result for projective schemes over surfaces.
\newblock  J. Reine Angew. Math. {\bf 606} (2007) 179-192.

\bibitem{CH}
S. Dale Cutkosky, J\"urgen Herzog.
\newblock Failure of tameness for local cohomology.
\newblock J. Pure Appl. Algebra {\bf 211} (2007), no. 2, 428-432. 

\bibitem{DV}
Mohammad T. Dibaei, Alireza Vahidi.
\newblock  Artinian and Non-Artinian Local Cohomology Modules.
\newblock  Canad. Math. Bull. {\bf 54} (2011), 619-629.


\bibitem{DNT}
K.Divaani-Aazar, R.Naghipour, M.Tousi.
\newblock  Cohomological dimension of certain algebraic varieties, 
\newblock  Proc. Amer. Math. Soc. {\bf 130} (2002), 3537-3544.

\bibitem{Go}
Roger Godement.
\newblock  Topologie alg\'ebrique et th\'eorie des faisceaux. 
\newblock  Actualit\'es Scientifiques et Industrielles, No. 1252. Hermann, Paris, 1973.

\bibitem{EGAIII}
Alexandre Grothendieck. 
\newblock  \'El\'ements de g\'eom\'etrie alg\'ebrique. III. \'Etude cohomologique des faisceaux coh\'erents. I. 
\newblock  Inst. Hautes \'Etudes Sci. Publ. Math. No. 11, 1961.

\bibitem{Jo}
Jean-Pierre Jouanolou.
\newblock  {\em Mono\"\i des}.
\newblock Publications de l'Institut de recherche math\'ematique avanc\'ee, Universit\'e Louis Pasteur, 1986.

\bibitem{Ke}
George R. Kempf.
\newblock  {\em Some elementary proofs of basic theorems in the cohomology of quasicoherent sheaves}.
\newblock  Rocky Mountain J. Math. {\bf 10} (1980), no. 3, 637--645.

\bibitem{No}
D.~G. Northcott.
\newblock {\em Finite free resolutions}.
\newblock Cambridge University Press, Cambridge, 1976.
\newblock Cambridge Tracts in Mathematics, No. 71.

\bibitem{Oo}
A. Ooishi.
\newblock Castelnuovo's regularity of graded rings and modules,
\newblock Hiroshima Math. J. {\bf 12} (1982), No. 3, 627-644.

\bibitem{Ra}
Ahad Rahimi.
\newblock Local cohomology of bigraded modules,
\newblock Ph. D. thesis, Universit\"at Essen, 2007.

\bibitem{St}
M. Stokes.
\newblock  Some dual homological results for modules over commutative rings, 
\newblock J. Pure Appl. Algebra {\bf 65} (1990), 153-162.

\bibitem{Su}
A. A. Suslin.
\newblock Triviality of certain cohomology groups,
\newblock  Zap. Nauchn. Sem. Leningrad. Otdel. Mat. Inst. Steklov. (LOMI) {\bf 94} (1979), 114--115.
\newblock {\it English translation :} Journal of Mathematical Sciences {\bf 19} (1982), No. 1, 1048-1049.

\end{thebibliography}
\end{document}